\documentclass[reqno, 11pt]{amsart}
\usepackage[utf8]{inputenc}
\usepackage[dvips]{graphicx,color}
\usepackage[all]{xy}
\usepackage{enumerate}
\usepackage[english]{babel}
\usepackage{graphicx,pstricks,pst-plot}
\usepackage{latexsym}
\usepackage{xcolor}
\usepackage[colorlinks=false, linkcolor=red, citecolor=green, urlcolor=blue]{hyperref}
\usepackage{setspace}
\usepackage{cancel}
\usepackage[left=2.5cm, right=2.5cm, bottom=2cm]{geometry}
\usepackage{amsmath, amssymb, amsthm, epsfig, mathtools}
\usepackage{hyperref, latexsym}
\usepackage[nameinlink,capitalise]{cleveref}
\usepackage{url} 
\usepackage{stmaryrd}
\usepackage{enumitem}

\usepackage{float}

\usepackage{tkz-base}
\usetikzlibrary{patterns} 
\usepackage{graphicx}
\usepackage{array}
\usepackage{tkz-fct}
\usepackage{geometry}
\usepackage{tikz}
\usetikzlibrary{shapes,arrows}
\usetikzlibrary{intersections}

\newcommand\bna{\begin{eqnarray*}}
\newcommand\ena{\end{eqnarray*}}
\newcommand\bnan{\begin{eqnarray}}
\newcommand\enan{\end{eqnarray}}

\newcommand{\Pc}{\mathcal{P}}
\newcommand{\Q}{\mathcal{Q}}
\renewcommand{\P}{\mathcal{P}}
\newcommand{\bC}{\mathbf{C}}

\newcommand{\al}{\alpha}

\newcommand{\veps}{\varepsilon}

\providecommand{\norm}[1]{\lVert#1\rVert}

\newcommand\nor[2]{\left\|#1\right\|_{#2}}

\usepackage{dsfont}
\usepackage{mathtools}

\def\today{\ifcase\month\or
  January\or February\or March\or April\or May\or June\or
  July\or August\or September\or October\or November\or December\fi
  \space\number\day, \number\year}
\DeclareMathOperator{\supp}{\mathrm{supp}}

\newtheorem{theorem}{Theorem}[section]
\newtheorem{lemma}[theorem]{Lemma}

\newtheorem{corollary}[theorem]{Corollary}
\theoremstyle{definition}

\theoremstyle{remark}

\newtheorem{assum}{Assumption} 
\newenvironment{assump}[2][]
  {\renewcommand\theassum{#2}\begin{assum}[#1]}
  {\end{assum}}

\newcommand{\mc}{\mathcal}

\newcommand{\M}{\Omega}

\newcommand{\R}{\mathbb{R}}
\newcommand{\N}{\mathbb{N}}
\newcommand{\Z}{\mathbb{Z}}

\newcommand{\intav}[1]{\mathchoice {\mathop{\vrule width 6pt height 3 pt depth  -2.5pt
\kern -8pt \intop}\nolimits_{\kern -6pt#1}} {\mathop{\vrule width
5pt height 3  pt depth -2.6pt \kern -6pt \intop}\nolimits_{#1}}
{\mathop{\vrule width 5pt height 3 pt depth -2.6pt \kern -6pt
\intop}\nolimits_{#1}} {\mathop{\vrule width 5pt height 3 pt depth
-2.6pt \kern -6pt \intop}\nolimits_{#1}}}

\begin{document}

\title[Propagation results for semilinear conservative PDE\lowercase{s}]{Global propagation of analyticity and unique continuation for semilinear conservative PDE\lowercase{s}}
\author[]{Camille Laurent}
\date{\today}
\address{CNRS UMR 9008, Universit\'e Reims-Champagne-Ardennes, Laboratoire de Math\'ematiques de Reims (LMR), Moulin de la Housse-BP 1039, 51687 REIMS cedex 2, France.}
\email{camille.laurent@univ-reims.fr}
\author[]{Crist\'obal Loyola}
\address{Sorbonne Université, Université Paris Cité, CNRS, Laboratoire Jacques-Louis Lions, LJLL, F-75005 Paris, France}
\email{cristobal.loyola@sorbonne-universite.fr}

\allowdisplaybreaks
\numberwithin{equation}{section}

\begin{abstract}
We review some recent results in which we develop a new method for proving global unique continuation for some conservative PDEs. The main tool is to prove some global propagation of analyticity. We first present some known results on the subject. Then, we sketch the abstract method we use, which relies on the property of finite determining modes. We give applications to semilinear wave, plates and Schr\"odinger equations. This note was written for the \emph{Proceedings of the Journées EDP 2025}.
\end{abstract}

\maketitle

\section{Introduction}

In this article, we intend to survey some recent results concerning some proofs of global propagation of analyticity that allow to obtain some global unique continuation results with natural global geometric assumptions.

In a first part we will present the problematic. Then, in Section \ref{s:known}, we will present some  known results on the subject. Afterward, in Section \ref{s:mainNLW} we will present our results for the semilinear wave equation. Then, in Section \ref{s:idea}, we will present an abstract result and give a sketch of the proof. Then we will present some further results concerning the semilinear plate equation in Section \ref{s:plates} and semilinear Schr\"odinger equation in Section \ref{s:NLS}.

The results  are sometimes written in a slightly informal way. We refer to the original articles \cite{LL:24, Loy25:NLS, Loy25:NLP} for more precise formulations.

\subsection{Propagation of information}

The informal question we address is the following : if we can say "something" about the behavior of a solution to a PDE, somewhere in a subdomain, what can we say about the solution on the whole domain?

From a physical perspective, this means that we have information only on part of the domain, and we ask whether it is possible to get the knowledge of our system everywhere. In more mathematical terms, this can be expressed through the following questions.

\medskip

Let $P$ be a linear or nonlinear differential operator defined  on $\Omega$ for a domain $\Omega\subset \R^n$ (or a manifold). Consider $u$ is a solution $P(u)=0$ in $\Omega$. For $\omega\subset \Omega$, we can study :
\begin{itemize}
\itemsep-0.25em
\item propagation of nullity (or uniqueness): unique continuation
\bna
 u|_{\omega}=0 \quad\overset{?}{\Longrightarrow} \quad u= 0\ \text{ on } \Omega .
\ena
\item propagation of smallness : quantification of the unique continuation
\bna
u\ \text{ small in }\omega \quad\overset{?}{\Longrightarrow} \quad u\  \text{ small in }\Omega.
\ena

\item propagation of regularity:
\bna
u\ \text{ regular in }\omega \quad\overset{?}{\Longrightarrow}\quad u\ \text{ regular in } \Omega.
\ena
\end{itemize}
The method we develop is quite general for conservative equations and relies on observability estimates. We provide an abstract result, presented in more detail in Section \ref{s:idea}, that allow us to propagate analyticity in time from the observation to the full solution. This result, or a possible refinement, yields applications to nonlinear wave equations (see Section~\ref{s:mainNLW}), nonlinear plate equations (see Section~\ref{s:plates}) and nonlinear Schr\"odinger equations (see Section~\ref{s:NLS}). We believe that it could be applied to several other systems and is amenable to generalizations.

Before getting into the details of the results, we discuss some motivations for the questions we address. Propagation of analyticity and unique continuation for nonlinear equations, beyond their intrinsic interest, arise mainly in various rigidity problems. In certain situations, we need to identify some objects that satisfy very specific properties and they are often expected to be the asymptotic solutions at large time. We present some examples below.
\begin{itemize}
    \item \emph{Nonlinear stabilization problems:} consider some partially damped equation of the form, for instance, $\partial_t^2u-\Delta u+\mathds{1}_{\omega}\partial_t u=f(u)$ or $\partial_t^2u+\Delta^2 u+\mathds{1}_{\omega}\partial_t u=f(u)$ with $\omega \subset \M$. A main obstacle to the decay would be the existence of undamped solutions that would satisfy $\partial_t u=0$ on $\omega$. It is very desirable to obtain such a result with some geometric assumptions on $\omega$ that are as minimal as possible. In the case that these objects do exist, they are expected to describe the asymptotic behavior of solutions, for instance if they belong to some attractor set. We refer, for instance, to \cite{Z:91,DLZ03,JL13,L:11} for such problems in the context of wave equations.
    \item \emph{The rigidity of stationary black holes:} this problem consists in classifying the set of solutions of the Einstein equation that contain a Killing vector field, are asymptotic to the Minkowski metric, and satisfy suitable additional assumptions. It is conjectured that there are no solutions of this type beyond the Kerr family. This subject has recently been the subject of intense research (see \cite{IK:15} for a survey). The result is already known to be true for analytic solutions. In particular, a result of propagation of analyticity for such solutions would certainly allow some great progress on this question. Although our paper is restricted to semilinear equations, we might expect that our result, along with additional ideas, could extend to the fully nonlinear setting.
    \item \emph{Nonlinear scattering theory:} on $\R^d\setminus \mathcal{O}$, where $\mathcal{O}$ is a finite union of compact obstacles, consider, for instance, equations $\partial_t^2u-\Delta u=f(u)$ or $\partial_t^2u+\Delta^2 u=f(u)$. The question of scattering asks whether there exists a linear solution  $u_L$ such that $u(t)-u_L(t)\to 0$ as $|t|\to+\infty$ in some suitable norm. There might be several obstacles to this property. One is the existence of solutions with the compactness property, which means that $\{u(t)\ |\ t\in\R\}$ is compact. An even worse scenario is the existence of a solution that remains supported in a fixed compact set in space. This is the typical situation of unique continuation that our analysis might cover. Note that the current methods for this problem use Morawetz-type estimates (see, for instance, \cite{BSS:09} outside of a star-shaped obstacle or \cite{KVZ:16} for NLS outside of a convex obstacle), but it seems that they cannot reach the more general geometries we wish to treat.
\end{itemize}
In order to be more concrete, let us present more precisely how the question of propagation of information reads in the context of the semilinear wave equation.

\subsection{Propagation of information for semilinear wave equation}
In many evolution equations, time has a special role and it is very useful to consider the questions we described before in this context. The equation under consideration is the following semilinear wave equation
\begin{align}\label{NLW}\tag{NLW}
    \left\{\begin{array}{rl}
        \partial_t^2 u-\Delta_g u+f(u)=0  &\ (t, x)\in [0, T]\times \Omega,\\
        u_{|_{ [0, T]\times\partial\Omega}}=0            &\ (t, x)\in [0, T]\times\partial\Omega,\\
        (u, \partial_t u)(0)=(u_0, u_1) &\ x\in \Omega,
    \end{array}\right.
\end{align}
and we aim to study the three following properties:
\begin{itemize}
\itemsep-0.25em
\item propagation of nullity (or uniqueness): unique continuation
\bna
 \partial_{t}u|_{[0,T]\times\omega}=0 \quad\Longrightarrow \quad u= 0  .
\ena
\item propagation of smallness : observability estimates
\bna
 \partial_{t}u\ \text{ small in }[0,T]\times\omega \quad\Longrightarrow \quad (u_0, u_1)\  \text{ small}.
 \ena

\item propagation of  analyticity:
\bna
u\ \text{ analytic (in time) in }[0,T]\times\omega \quad \Longrightarrow \quad u\ \text{ analytic in } [0,T]\times\Omega.
\ena
\end{itemize}
The main result we present concerning semilinear wave equations will be to provide a positive answer to these questions when $(\omega,T)$ satisfy the Geometric Control Condition, see Theorem \ref{thm:unique-continuation-nlw} and \ref{thm:analytic-prop} below.

Before addressing these results in detail, let us discuss the known results on unique continuation and explain the difficulties in obtaining global results under natural geometric assumptions.

\section{Known results of unique continuation}
\label{s:known}

\subsection{Classical local unique continuation theorems for linear equations}
\label{s:classic}
Many of the available unique continuation results are local in nature. We consider $\Omega$ a bounded open subset of $\R^n$, $P$ a linear differential operator on $\Omega$, $x_0 \in \Omega$ a point, and a hypersurface $S=\left\{\Phi=0\right\}$ containing $x_0$. The aim is to prove some local unique continuation for an operator $P$ across the hypersurface $S=\left\{\Phi=0\right\}$, say, a statement like
\begin{equation*}
Pu = 0  \text{ in } \Omega , \quad u  = 0  \text{ in }\Omega \cap \left\{\Phi\geq 0\right\}
\overset{?}{\Longrightarrow}  u=0\text{ close to }x_0.
\end{equation*}
In particular, we want to prevent the situation in which a smooth function $u$ both solves $Pu=0$ and vanishes on $S$, possibly ``flatly'', in the sense that all its derivatives vanish. This problem has a very long history, and it would be impossible to be fair with the full literature. We refer for instance to \cite{LernerBook19,LL:book} for a general treatment or \cite{LL:23notes} for presentation with special emphasis on the wave equation. We can roughly summarize the types of results in the following Figure \ref{tableUCP}, where $p$ denotes the principal symbol of the linear operator $P$.

\begin{figure}
  \begin{tabular}{|m{2.2cm}||m{3.5cm}|m{4cm}|m{3.5cm}|}
\hline
Theorem&Holmgren~\cite{Holmgren}, John~\cite{John:49}&Tataru~\cite{Tat95,Tat99}, Robbiano-Zuily~\cite{1998:robbiano-zuily:ucp-analytic}, H\"ormander~\cite{Hor:97}& Carleman \cite{C:39}, Calder\'on~\cite{C:58}, H\"ormander~\cite{H:63}\tabularnewline
   \hline
 Regularity of coefficients & analytic coefficients & partially analytic in one variable $x_a$& coefficients $C^{\infty}$ (or even $C^1$) \tabularnewline
    \hline
Geometric assumptions on the hypersurface&$\Phi$ non characteristic for $P$: $p(x,\nabla \Phi)\neq 0$&$\Phi$ strictly pseudoconvex in $\{\xi_a = 0\}$ (OK for waves and $\Phi$ non characteristic)&  $\Phi$ strictly pseudoconvex, see \eqref{pseudoconvexe} \tabularnewline
\hline

\end{tabular} 
\caption{The main local unique continuation theorems}
\label{tableUCP}
\end{figure}
\bigskip
The oldest results on local unique continuation are of Holmgren-John type (Column 1 in Figure \ref{tableUCP}). They ensure the local unique continuation for linear operators with analytic coefficients under the assumption that $\Phi$ is non characteristic, that is 
\begin{equation}
\label{noncaract}
    p(x,\nabla \Phi(x))\neq 0.
\end{equation} For the specific example of the wave operator $P=\partial_t^2-\Delta$ with principal symbol $p(x,\xi)=-\xi_t^2+|\xi_x|^2$, where $\xi_t$ (resp. $\xi_x$) are the Fourier variable dual to the time variable $t$  (resp. the space variable $x$), it can be written
\begin{equation*}
(\partial_t\Phi)^2\neq |\nabla_x \Phi|^2 \ \Longleftrightarrow\ \{\Phi=0\} \text{ not tangent to the light cone.}
\end{equation*}
This last assumption can be seen to be almost optimal for the 1D wave equation by taking the example of a solution with initial datum supported in $\{x, |x|\leq 1\}$ and with support the truncated light cone $\{(t,x), |x|\leq |t|+1\}$. In that case, the unique continuation along the hypersurface $\{x=t+1\}$ is false. Yet, the assumption of analytic coefficients is very demanding and many authors tried to go further. 

\bigskip

One of the main tools in this context is the use of Carleman estimates (Column 3 in Figure \ref{tableUCP}). It has led to a general theory (see for instance \cite[Section 18]{Hoermander:V4} or \cite{Lerner:19}) where the crucial assumption for the hypersurface is that  $\Phi$ strictly pseudoconvex. It is a quite complicated assumption, which, in the case of differential operator of order $2$ with real principal symbol slightly simplifies to
\begin{equation}
 \label{pseudoconvexe}   
p(x,\xi)=\{ p ,\Phi\}(x, \xi)=0 \quad \Longrightarrow \quad \{p,\{p,\Phi\}\}(x, \xi)>0  \quad \text{for all } \xi \neq 0,
\end{equation}
where $\{,.\}$ is the Poisson bracket. Roughly speaking it is a condition of order $2$ that ensures that 
the bicharacteristic of $p$ that are tangent to $\{\Phi= 0\}$ stay on the side $\{\Phi\geq 0\}$. 
 For instance, for the wave equation, the pseudoconvexity condition of a hypersurface $\{\Phi=0\}$ for the wave operator writes:
\begin{align*}
X_{t}^{2}=|X_{x}|^{2}\text{ and }d\Phi(x_0)(X)=0  \ \  \Longrightarrow \ \  \text{Hess}~ \Phi(x_0) (X , X) >0  \ \text{ for all } X=(X_{t},X_{x}) \in  \R^{1+d} \setminus \{0\}.
\end{align*}
For instance, in the particular case $\Phi(t,x)=\Phi(x)$, that is the unique continuation from a cylinder, it is sufficient to ask the (usual) convexity assumption $\text{Hess}~\Phi>0$. If $\Phi(x)=|x|^2-1$, $\Phi$ works for unique continuation, but not $-\Phi$. That means that you can prove unique continuation for the wave equation from the exterior of a ball to the interior, but not from inside to outside. This is in strong contrast with the Holmgren-John theorem which did not have preference directions.

It turns out that some geometric assumptions stronger than the Holmgren theorem are actually needed in general. Indeed, concerning the unique continuation with even smooth $V$, the classical counterexamples of Alinhac-Bahouendi~\cite{AB:95} (see also the earlier \cite[Theorem
8.9.4]{H:63}), refined by H\"ormander \cite{H:00}, are quite striking. They show that H\"ormander's pseudoconvexity condition is not far from being optimal for local unique continuation. For any $s>1$ and $d\geq 2$, they construct some $u$ and $V \in \mathcal{G}^s(B_{\R^{1+d}}(0,1), \mathbb{C})$ (Gevrey functions) so that 
\begin{align*}
\partial_{t}^{2}u-\Delta u&=V u \textnormal{ on }B_{\R^{1+d}}(0,1),\\
\supp(u)&=\left\{(t,x_{1},\dots, x_{d})\ |\  x_{1}\geq 0\right\}\cap B_{\R^{1+d}}(0,1).
\end{align*}
This suggests that in geometrical situations where the strong pseudoconvexity of the hypersurface is not satisfied, we cannot expect local unique continuation for a potential $V$ that is not analytic.

\medskip

The last type of results concern the unique continuation with partial analyticity  (Column 2 in Figure \ref{tableUCP}) which will be crucial for obtaining our main result Theorem \ref{thm:unique-continuation-nlw} from the propagation of regularity result Theorem \ref{thm:analytic-prop}. The history of this theory is quite long, with several breakthroughs and improvements, in particular by Tataru~\cite{Tat95,Tat99}, Robbiano-Zuily~\cite{1998:robbiano-zuily:ucp-analytic} and H\"ormander~\cite{Hor:97}, that we do not detail here (see for instance \cite[Section 4]{LL:23notes} for more details).

In this context, $\R^n$ is decomposed into $\R^{n_a}\times \R^{n_b}$ and we assume that the coefficients of the operator $P$ are analytic in the variable $x_a$. In most of our examples, we will use $n_a=1$ and $x_a=t$ will be the time variable. There is now a general theory comparable to the one for classical Carleman estimates which involves some related assumptions of strictly pseudoconvex functions in the cone $\{\xi_a=0\}$ which roughly requires to check the usual pseudoconvexity assumption (which, in general, contain \eqref{pseudoconvexe} and other convexity type requirements) only for the $\xi$ so that $\xi_a=0$, that is of the form $(0,\xi_b)$. We refer for instance to \cite{Lerner:19} or \cite{LL:book} for a more general overview. It turns out that for the wave equation, this assumption reduces to the noncharacteristic assumption \eqref{noncaract}, that is $p(x,\nabla \Phi)\neq 0$, see for instance \cite[Lemma 3.6]{LL:23notes} for a proof. So, as already mentioned for the Holmgren-John theorem, concerning the geometry constraints, this result is very satisfactory. It remains the assumption of coefficients analytic with time, but it is already a weaker assumption that the analytic assumption of the Holmgren-John result.

Also, it is worth mentioning that this theory actually encompasses the two types of results described before: the Holmgren-John's theorem is a particular case of this theory with $n_{a}=n$ while the H\"ormander's theorem can be deduced from the case $n_{a}=0$. It is actually interpolating between Holmgren-John's theorem and the H\"ormander's theorem.

\medskip

Finally, we might expect to prove directly a unique continuation theorem for nonlinear equation. Nevertheless, another counterexample that undermines such approach is provided by M\'etivier in \cite{M:93}, where it is proved that a nonlinear version of the Holmgren-John theorem fails in general. The operators for which it applies are not wave operators, but a nonlinear Holmgren-John theorem, even for a more specific class of operators has, up to our knowledge, never been obtained so far, except for scalar operators of order $1$.

\subsection{Some previous global results for semilinear waves}\label{s:UCPwaveglobal}
All the results presented in the previous subsection are local in nature. But it is of course possible to make them global by iterating them and constructing some foliation of hypersurface to finally obtain some global unique continuation theorem. Yet, the global geometry enters into the problem. We describe some of these issues and results for the nonlinear wave equation \eqref{NLW}.

From now onward, unless stated otherwise, $\Omega$ is a smooth compact Riemannian manifold with (or without) boundary of dimension $d$.
\subsubsection{Nonlinear solutions are also linear solutions}
\label{s:NLisL}
As already mentioned, a nonlinear version of the Holmgren-John theorem is known to be false in general, see~\cite{M:93}. Indeed, there are very few unique continuation results that genuinely take the nonlinearity itself into account. We are only aware of Alexakis-Shao~\cite{AS:15} that use the sign of the nonlinearity in a Carleman estimates (see also \cite{LBM:25} for the nonlinear Sch\"odinger equation).

In most of the articles, the strategy is the following "trick" to consider a nonlinear equation as a linear one.

Let $u$ be a solution of \eqref{NLW} and satisfying $ \partial_{t}u|_{[0,T]\times\omega}=0$. We want to know if it implies
\bnan
\label{UCP}
u= 0 \text{ on } [0,T]\times\Omega.
\enan
The idea for this problem is to consider the function  $w=\partial_{t}u$. It is solution of 
\begin{align*}
    \Box w+Vw&=0,\\
    w|_{[0,T]\times\omega}&=0,
\end{align*}
 with $V=f'(u)$. We can therefore see $w$ as a solution of $P_{V}w=0$ where $P_{V}$ is the linear operator $\Box +V$.
 But there is of course a price to pay. Even if the the coefficients of the nonlinear equation are very regular, the regularity of the potential $V=f'(u)$ now depends on the regularity of the solution itself. Therefore, in the classical unique continuation results presented in Section \ref{s:classic}, we need to inspect the allowed regularity of the lower order terms of the operator $P$ to make them coincide with the regularity of $V=f'(u)$.

Note that this "trick" is quite specific to the context of unique continuation where we don't look for constructing a nonlinear solution since is already given in the problem.
\subsubsection{ Using H\"ormander theorem and classical Carleman estimates}\label{s:carlemanhormander}
As already mentioned, the classical Carleman estimates and the H\"ormander theorem have the advantage that the regularity required is not so demanding. For instance, the usual theory allows some lower order terms that are $L^{\infty}$ and some additional dispersive type estimates, similar to Strichartz estimates, allow sometimes to lower this to $L^p([0,T],L^q(\Omega)$ spaces with suitable $(p,q)$. The limitation comes more from the geometric assumptions which are not  so natural. Basically, we need to be in some situations where there exists a foliation of hypersurfaces that are pseudoconvex as in \eqref{pseudoconvexe}. Achieving global unique continuation would then require iterating some local unique continuation results across a hypersurface $\{\Phi=0\}$. Yet, for a general configuration, the global geometric assumptions resulting from the use of  H\"ormander theorem for the unique continuation are not very natural and are stronger than \eqref{assumGCC}. Typical geometric assumptions are often of ''multiplier type'' (or Morawetz type), meaning $\omega$ is a neighborhood of $\left\{x\in \partial\Omega \left|(x-x_0)\cdot n(x)>0\right.\right\}$ which are known to be stronger than the geometric control condition (see \cite{M:03} for a discussion about the links between these assumptions). Moreover, on curved spaces, this type of condition often needs to be checked by hand in each situation, which is mostly impossible in general.

In regards to classical Carleman estimates, many authors explored the global assumptions needed to obtain them, as well as their consequences. Global Carleman estimates for waves were proved in~\cite[Chapter~4]{FI:96} and \cite{BDBE:13} with applications to controllability and inverse problems. Another line of investigations in a geometric context was undertaken in~\cite{DZZ:08,Shao:19,JS:21} aiming to present geometric assumptions that would ensure the usual pseudoconvexity, with applications to observability estimates and null-controllability.

For treating nonlinear problems, admitting lower-order terms in unique continuation results is crucial. For instance, in the present work, a nonlinearity of the form $f(u)$ is treated as a term $V u$ with potential $V$. Here, a previous analysis allows to obtain that the nonlinear solution is actually very regular. Yet, having nonlinear problems in mind, it has sometimes been a goal to minimize the regularity of the admissible lower-order terms. A global unique continuation statement was proved in~\cite{Ruiz:92}, with application to energy decay for nonlinear waves. This led to some ``dispersive'' Carleman estimates with Strichartz-type spaces. The literature is vast, and we refer, for instance, to~\cite{KRS:87,DDSF:05, KT:05}, and the references therein. For a more complete historical overview of Carleman estimates, we refer to \cite[Section 4]{LL:23notes}.

\subsubsection{Using unique continuation with partial analyticity}
\label{s:UCPpartial}
To globalize the unique continuation for wave equations from the Tataru-Robbiano-Zuily-H\"ormander theory, we only need to construct a foliation of noncharacteristic hypersurfaces, which is much simpler than constructing pseudoconvex hypersurfaces. Such foliations can be found, for instance, in \cite{L:92} or \cite{L:00}. This approach yields the following global result.
\begin{theorem}[Tataru-Robbiano-Zuily-H\"ormander]\label{thm:qucp}
    Let $V$ be bounded and analytic in $t$. Let $\omega\subset \Omega$ open and $T>2\sup_{x\in\Omega}\text{dist}(x,\omega )$. Let $(u_0, u_1)\in H_0^1\times L^2(\Omega)$ and $u$ solution of
    \begin{align*}
        \left\{\begin{array}{ll}
\partial_{t}^{2}u-\Delta_{g} u+Vu=0&\ \text{ in }\ (0, T)\times\Omega\\
            u_{|_{\partial\Omega}}=0 &\ \text{ in }\ (0, T)\times \partial\Omega\\
            (u, \partial_t u)(0)=(u_0, u_1), &\ 
        \end{array} \right.
    \end{align*}
that satisfies $u=0$ on $[0,T]\times \omega$. Then $u=0$ on $[0,T]\times \Omega$.
\end{theorem}
This theorem is very satisfactory from a geometric point of view, since no assumption is made on the observation set $\omega$. Moreover, the limit time $2\sup_{x\in\Omega}\text{dist}(x,\omega )$ can be seen to be optimal by finite speed of propagation: we can consider initial data compactly supported close to a point $x_1$ so that $\text{dist}(x,\omega )>\sup_{x\in\Omega}\text{dist}(x,\omega )-\varepsilon$ and work on the interval $[-T,T]$.

Note also that a quantitative version of Theorem \ref{thm:qucp} has also been given in Laurent-L\'eautaud~\cite{LL:19} giving some estimates of the form 
\bna
\nor{(u_0,u_1)}{H^1\times L^2}\leq C e^{C\Lambda}\nor{u}{L^2((0,T); H^1(\omega))},
\ena
where $\Lambda=\frac{\nor{(u_0,u_1)}{H^1\times L^2}}{\nor{(u_0,u_1)}{L^2\times H^{-1}}}$ needs to be thought as the typical frequency of the initial datum.

Returning to our unique continuation problem for \eqref{NLW}, and following the discussion of Section \ref{s:NLisL}, the main obstruction lies in the fact that $V=f'(u)$ must be analytic in time so that Theorem \ref{thm:qucp} can be applied. But, unlike parabolic equations, solutions to hyperbolic equations are in general far from being analytic. This explains why the fundamental step of our proof is to prove that such analyticity actually holds for the very particular solutions satisfying \eqref{UCP}. This is the content of Theorem \ref{thm:analytic-prop} below. But before getting tho this result, we need to introduce the observability property, which is a quantitative version of unique continuation.

\subsection{Observability under Geometric Control Condition}
An observability estimate is, in some sense, the strongest possible quantification of the unique continuation. It writes as a linear estimates where the norm of the observation is bounded below by the norm of the initial data. When it turns to the linear wave equation, the right geometric assumption is the following classical Geometric Control Condition which is very natural since we learn in physics that light travels along the rays of geometric optics.

\begin{assump}{GCC}\label{assumGCC}
 We say $(\omega,T)$ Geometric Control Condition \eqref{assumGCC} if  every ray of geometric optics traveling at speed $1$ meets $\omega$ in time $<T$.
\end{assump}

The definition of "ray of geometric optics" is actually not trivial to define mathematically, see Melrose-Sj\"ostrand~\cite{MS:78}, but quite natural to draw, see Figure \ref{fig:GCC} below. It leads to the following observability estimate.
\begin{theorem}[Bardos-Lebeau-Rauch~\cite{BLR92}]\label{t:BLR}
If $(\omega,T)$ satisfies \eqref{assumGCC}, then for every solution of
    \begin{align*}
        \left\{\begin{array}{ll}
\partial_{t}^{2}u-\Delta_{g} u=0&\ \text{ in }\ (0, T)\times\Omega\\
            u_{|_{\partial\Omega}}=0 &\ \text{ in }\ (0, T)\times \partial\Omega\\
            (u, \partial_t u)(0)=(u_0, u_1), &\ 
        \end{array} \right.
    \end{align*}
we have
\bna
\nor{(u_0,u_1)}{H^1\times L^2}\leq C \nor{\partial_{t}u}{L^2([0,T]\times\omega)}.
\ena
\end{theorem}

\begin{figure}[H]
\label{fig:GCC}
\tikzset{every picture/.style={line width=0.75pt}}

\begin{tikzpicture}[scale=1, x=0.75pt,y=0.75pt,yscale=-1,xscale=1]
\useasboundingbox (90,15) rectangle (470,245);

\draw   (357.94,94.05) .. controls (383.54,94.05) and (392.75,87.05) .. (411.46,85.64) .. controls (430.16,84.24) and (445.23,87.05) .. (459.26,102.96) .. controls (473.3,118.87) and (472.78,153.49) .. (450.95,168.47) .. controls (429.12,183.44) and (404.33,182.34) .. (364.34,192.57) .. controls (324.35,202.81) and (292.84,224.08) .. (272.58,233.46) .. controls (252.32,242.85) and (220.38,246.32) .. (199.59,219.45) .. controls (178.79,192.57) and (162.54,190.56) .. (152.15,185.5) .. controls (141.75,180.45) and (95.62,160.58) .. (92.42,108.12) .. controls (89.22,55.66) and (118.36,35.87) .. (145.65,31.78) .. controls (172.93,27.69) and (175.96,56.35) .. (194.15,55.57) .. controls (212.34,54.79) and (215.16,30.23) .. (238.75,22.81) .. controls (262.35,15.39) and (278.2,31.38) .. (302.2,55.7) .. controls (326.19,80.01) and (332.35,94.05) .. (357.94,94.05) -- cycle ;
\draw [color={rgb, 255:red, 74; green, 134; blue, 226 }  ,draw opacity=1 ]   (184.07,196.78) .. controls (234.19,293.18) and (317.37,201.9) .. (333.36,189.96) ;
\draw [color={rgb, 255:red, 74; green, 134; blue, 226 }  ,draw opacity=1 ]   (367.74,160.09) -- (333.36,189.96) ;
\draw [color={rgb, 255:red, 74; green, 134; blue, 226 }  ,draw opacity=1 ]   (103.03,58.59) -- (184.07,196.78) ;
\draw [color={rgb, 255:red, 74; green, 134; blue, 226 }  ,draw opacity=1 ]   (103.03,58.59) -- (298.17,51.76) ;
\draw [color={rgb, 255:red, 74; green, 134; blue, 226 }  ,draw opacity=1 ]   (298.17,51.76) -- (263.96,118.32) ;
\draw  [color={rgb, 255:red, 74; green, 134; blue, 226 }  ,draw opacity=1 ] (271.3,114.55) -- (263.72,118.36) -- (262.29,110.65) ;
\draw [color={rgb, 255:red, 74; green, 74; blue, 74 }  ,draw opacity=1 ] [dash pattern={on 0.84pt off 2.51pt}]  (275.25,100.08) .. controls (280.88,101.45) and (310.02,105.76) .. (331.78,99.6) .. controls (353.55,93.45) and (372.71,80.25) .. (375.31,66.8) ;
\draw [color={rgb, 255:red, 192; green, 134; blue, 233 }  ,draw opacity=1 ][fill={rgb, 255:red, 192; green, 134; blue, 233 }  ,fill opacity=0.93 ]   (143.13,68.74) .. controls (176.91,53.24) and (179.51,106.17) .. (217.19,134.25) .. controls (254.86,162.33) and (265.26,154.14) .. (296.44,138.93) .. controls (327.62,123.72) and (404.27,95.65) .. (428.95,103.83) .. controls (453.63,112.02) and (431.55,159.99) .. (417.26,158.82) .. controls (402.97,157.65) and (425.32,129.33) .. (412.06,123) .. controls (398.8,116.67) and (297.74,181.98) .. (270.45,185.72) .. controls (243.17,189.47) and (239.27,186.62) .. (215.89,177.98) .. controls (192.5,169.35) and (108.06,87.16) .. (143.13,68.74) -- cycle ;

\draw (99.39,16.9) node [anchor=north west][inner sep=0.75pt]   [align=left] {$\displaystyle \Omega $};

\draw (312.26,17.63) node [anchor=north west][inner sep=0.75pt]  [font=\small] [align=left] {\begin{minipage}[lt]{173.97pt}\setlength\topsep{0pt}
\begin{center}
{\fontfamily{pcr}\selectfont Projection of a generalized bicharacteristic}\\{\fontfamily{pcr}\selectfont onto the variable }$\displaystyle x$
\end{center}

\end{minipage}};

\draw (244.6,201.73) node [anchor=north west][inner sep=0.75pt]    {$\textcolor[rgb]{0.56,0.07,1}{\omega }$};

\end{tikzpicture}
\caption{The Geometric Control Condition of Bardos-Lebeau-Rauch.}
\end{figure}
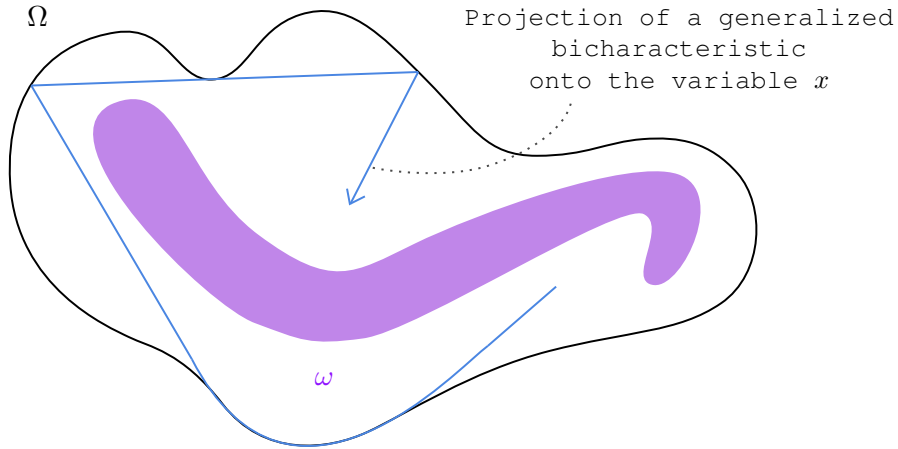

\section{Main results for semilinear waves}
\label{s:mainNLW}
\subsection{Unique continuation for semilinear waves}
When it comes to the semilinear equation, our main result is the following.
\clearpage
\begin{theorem}[Unique continuation, Laurent-Loyola~\cite{LL:24}]\label{thm:unique-continuation-nlw}
$d\leq 3$. Assume that:
    \begin{itemize}
\item $(\omega, T)$ satisfies the Geometric Control Condition \eqref{assumGCC},
    \item $f$ is energy subcritical and real analytic, $f(0)=0$.
    \end{itemize}
     If $(u,\partial_t u)\in C^0([0, T], H^1_0\times L^2(\M))$ solution  of \eqref{NLW} with finite Strichartz norms satisfies {$\partial_t u=0$} in $[0, T]\times \omega$, then $\partial_t u=0$ in $[0, T]\times \Omega$ and $u$ is a stationary solution (solution of elliptic equation).  If, moreover, the nonlinearity is defocusing
    \begin{align}
    \label{defdefocusing}
        sf(s)\geq 0\quad  \textnormal{ if }\partial \M\neq \emptyset,\\
    \nonumber    sf(s)\geq \gamma s^2\quad \textnormal{ if }\partial \Omega= \emptyset,
    \end{align}
    for some $\gamma>0$ and for all $s\in \R$, then $u\equiv 0$.
\end{theorem}

It was already proved in infinite time $T=+\infty$ by Joly-Laurent~\cite{JL13} using an infinite time regularization argument due to Hale-Raugel~\cite{HR03}. It was also proved in infinite time in the presence of weak trapping by Joly-Laurent~\cite{2020:joly-laurent:decay-nlw-no-gcc}, under more restrictive assumptions on the nonlinearity. The novelty here is the fact that it holds in the exact time of Geometric Control Condition. Moreover, the new method allows to prove some propagation of analyticity and is quite flexible to treat other models; see Sections \ref{s:plates} and \ref{s:NLS}. Once this unique continuation is established, it is possible to derive observability estimates for the \textbf{nonlinear} equation \eqref{NLW},
\bna
\nor{(u_0,u_1)}{H^1\times L^2}\leq C \nor{\partial_{t}u}{L^2([0,T]\times\omega)},
\ena
for the subcritical defocusing wave equation ($C$ depends on the size of the data).

As explained earlier in Section \ref{s:UCPwaveglobal}, the main step consists in proving that the solutions considered in Theorem \ref{thm:unique-continuation-nlw} are actually analytic in time. This property allows us to apply the general theory of unique continuation under partial analyticity and conclude.
\subsection{Propagation of analyticity}
Our second main result establishes that a solution of \eqref{NLW} which is analytic in time on a subdomain $[0,T]\times \omega$ satisfying the Geometric Control Condition is actually analytic in time on the whole domain. More precisely, we have the following.
\begin{theorem}[Propagation of analyticity, Laurent-Loyola~\cite{LL:24}]\label{thm:analytic-prop}
    Let $d\leq 3$ and $\sigma\in (1/2, 1]$. Let $(u,\partial_t u)\in C^0([0, T], H^{1+\sigma}\cap H^1_0\times H^{\sigma}_0(\M))$ be a solution of \eqref{NLW}. Assume:
    \begin{enumerate}
        \item $(\omega, T)$ satisfies the Geometric Control Condition~\eqref{assumGCC}.
        \item $u$ analytic in $t$ on $(0,T)\times \omega$: the map $t\in (0, T)\mapsto \chi u(t, \cdot)\in H^{1+\sigma}(\M)$ is analytic for any $\chi\in C_c^\infty(\M)$ whose support is contained in $\omega$.
        \item $s\mapsto f( s)$ is real analytic and $f(0)=0$.
    \end{enumerate}
    Then $t\mapsto u(t, \cdot)$ is analytic from $(0, T)$ to $H^{2}(\M)\cap H_0^1(\M)$.
\end{theorem}
The proof of this theorem is actually a consequence of an abstract result described in Section \ref{s:abstract}. Note that, in order to apply Theorem \ref{thm:analytic-prop} to the solutions considered in Theorem \ref{thm:unique-continuation-nlw}, it is necessary to gain regularity from $H_0^1\times L^2$ to $H^{1+\sigma}\cap H_0^1\times H^{\sigma}_0(\M)$. This is possible using the Geometric Control Condition~\eqref{assumGCC} and methods inspired by Dehman-Lebeau-Zuazua~\cite{DLZ03}.

The previous theorem concerns analyticity in time only. Exploiting the hyperbolic property of the wave, it is possible to propagate analyticity to all variables.
\begin{corollary}
\label{coranalytic}
If, moreover, the metric and the boundary are analytic in $x$, then $u$ is analytic in all variables.
\end{corollary}
The proof of Corollary \ref{coranalytic} from Theorem \ref{thm:analytic-prop} is actually quite short. 

Since $u$ is analytic in time, we can define its extension in a small complex neighborhood of the time variable. For a fixed $t_0$, consider $v(s,x)=u(t_0+is,x)$. It is solution of the elliptic equation $-\partial_s^2 v-\Delta_g v+f(v)=0$. Since all the coefficients are analytic, $v$ is therefore analytic by some classical results for elliptic equations \cite{Friedman:58}.
\subsection{(ultra)Short bibliography on propagation of regularity}
The literature for the propagation of regularity and local results is huge. So, we only intend to give the lines of results that exist. \medskip

The problem is better understood in the $C^{\infty}$ or $H^s$ context. For the linear equation, the propagation of $H^s$ regularity along the rays of Geometric Optics for the Dirichlet problem is well understood since the work of Melrose-Sj\"ostrand \cite{MS:78}; see \cite{H:07} for a complete historical overview on the internal and boundary problem. Concerning the $H^s$ regularity of nonlinear equations, the microlocal propagation has been the object of several studies since the work of Bony \cite{B:81}. The global propagation from a set satisfying \eqref{assumGCC} is proved in Dehman-Lebeau-Zuazua~\cite{DLZ03} using a bootstrap argument and propagation of $H^s$ regularity for smoother source terms.

Regarding the propagation of analyticity for linear equation, the propagation along bicharacteristics was proved in  H\"ormander~\cite{H:71} for a large class of operators. Yet, the problem becomes more complicated in the presence of boundary, especially for glancing rays, where there can be some diffraction, see  Friedlander-Melrose~\cite{FM:77} or Rauch-Sj\"ostrand~\cite{RS:81} for details. The propagation of analyticity that we prove may also be of interest in this context.

For the propagation of analyticity for nonlinear equations, there are much less results dating back to the 1980s and 1990s. Alinhac-M\'etivier proved in \cite{AM:84,AM:84b} that if $u$ is a regular enough solution of a general nonlinear PDE, the analyticity of $u$ propagates along any hypersurface for which the real characteristics of the linearized operator cross the hypersurface transversally. Subsequently, there has been some activity to understand what kind of singularities propagate for nonlinear waves. It was found that the situation is quite complicated since microlocal analytic singularities do not remain confined to bicharacteristics as in the linear case, but can give rise to nonlinear interactions. For more details, see Godin~\cite{G:86} and G\'erard~\cite{G:88}. Yet, in our geometric context, obtaining a global result from local propagation of singularities typically involves propagation from hypersurfaces of the form $S=\{\Phi=0\}$ with $\Phi(t,x)=\Phi(x)$. In such cases, the operator is never hyperbolic with respect to $S$, and there can be some bicharacteristics transverse to $S$ as soon as $d\geq 2$. In particular, the results from \cite{AM:84,AM:84b} do not seem to apply.
\section{Idea of proof}
\label{s:idea}
\subsection{Abstract result}
\label{s:abstract}
Our assumption will be as follows.
\begin{assump}{1}
\label{assum}
\begin{enumerate}
\item $A$ is skew-adjoint with compact resolvent, giving ''Sobolev spaces'' $X^{\sigma}$.
\item \label{condobser}observability for $\mathbf{C}\in\mathcal{L}(X^\sigma, X^\sigma)$: 
\begin{equation}
\label{observabs}
    \nor{W_0}{X^\sigma}^2\leq \mathfrak{C}_{\text{obs}}^2\int_0^T \nor{\mathbf{C} e^{tA}W_0}{X^\sigma}^2dt.
\end{equation}
\item \label{Fcomp}$F$ is a ''compact'' nonlinearity with a gain of derivatives : $X^{\sigma}\mapsto X^{\sigma+\varepsilon}$.
\item \label{Fanalytic} $F$ is analytic.
\item a technical assumption: there exists $s>0$ so that $[(A^*A)^s,\mathbf{C}]\in \mc{L}(X^{\sigma+2s}, X^{\sigma+\veps})$.
\end{enumerate}
\end{assump}

The abstract result we obtain is the following. For simplicity, we state the result that proves the propagation of the analytic regularity for solution whose observation vanishes, which is sufficient for the unique continuation property of Theorem \ref{thm:unique-continuation-nlw}. The result used to prove Theorem \ref{thm:analytic-prop} involves additional terms, but the spirit of the proof is the same. We believe that this simplified formulation makes the proof easier to present.
\begin{theorem}[Laurent-Loyola~\cite{LL:24}]\label{thmabstractanalyticintro}
Let $T^{*}>T$. 

Then, any solution $U\in C^{0}([0,T^{*}],X^{\sigma})$ satisfying 
\begin{align}
\label{e:UCPabstraite}
    \left\{\begin{array}{cr}
         \partial_t U=AU+F(U)& \textnormal{ on } [0,T^{*}], \\
         \mathbf{C} U(t)=0 &\textnormal{ for }t\in [0,T^{*}],
    \end{array}\right.
\end{align}
is real analytic in $t$ in $(0,T^*)$ with value in $X^{\sigma}$.
\end{theorem}
Let us now comment how to apply it to the semilinear wave equation.
Let $\beta\geq 0$. Set $X=H_0^1(\M)\times L^2(\M)$ and introduce the operator 
\begin{align*}
    A=\begin{pmatrix}
    0 & I \\
    \Delta_g-\beta & 0
    \end{pmatrix}\ \text{ with }\ D(A)=\big(H^2(\M)\cap H_0^1(\M)\big)\times H_0^1(\M).
\end{align*}
Here, $\Delta_g$ denotes the Dirichlet Laplacian defined on $L^2$
For the standard scalar product on $X$ (with $\norm{u}_{H^1_0}^2=\norm{\nabla u}_{L^2}^2+\beta \norm{u}_{L^2}^2$),
where we used the same notation for the Dirichlet Laplacian defined on $H^1_0(\M)$ and $L^2(\M)$ with their natural respective domains. For $\sigma\geq 0$, $H^\sigma_D$ and  $X^\sigma$ denote the spaces 
\begin{align*}
H^\sigma_D&=D((-\Delta_g)^{\sigma/2}),\\
 X^\sigma&= D((A^*A)^{\sigma/2})=  D((-\Delta_g)^{(1+\sigma)/2})\times D((-\Delta_g)^{\sigma/2})=H^{1+\sigma}_D\times H^{\sigma}_D.
\end{align*}

For $(\omega, T)$ satisfying \eqref{assumGCC}, the observation operator $\bC\in \mc{L}(X^\sigma, X^\sigma)$ is given by
\begin{align*}
    \bC(\phi, \psi)=(0, b_\omega\psi),
\end{align*} where $b_\omega$ is a smooth function so that $b_{\omega}(x)\geq 1$ for $x\in \omega$. Condition~\eqref{condobser} of Assumption~\ref{assum} is exactly Theorem~\ref{t:BLR}. In regards to the nonlinearity, we set
\begin{align*}
    F: (u, v)\in X^\sigma\longmapsto \big(0, f(u)\big)\in X^\sigma.
\end{align*}
One the one hand, we can roughly say that part \eqref{Fcomp} of Assumption \ref{assum} comes from the fact that $f(u)\in H_D^\sigma$ whilst $u\in H_D^{1+\sigma}$, so a small gain is possible while working at a fixed regularity scale. On the other hand, $F$ inherits the regularity of $f$, which verifies part \eqref{Fanalytic} of Assumption~\ref{assum} provided that $f$ is analytic. In this formulation, a solution to the semilinear wave equation \eqref{NLW} satisfying $\partial_t u=0$ on $[0, T]\times\omega$, can be written in the abstract form \eqref{e:UCPabstraite} with state $U=(u, \partial_t u)$.

\subsection{Scheme of proof}\label{s:scheme}
The proof of Theorem~\ref{thmabstractanalyticintro} is partially inspired by Hale-Raugel~\cite{HR03} in the context of the regularity of the attractor. The idea is to prove the property of finite determining modes: for this particular class of solutions, the knowledge of the low-frequency component of the solution allows us to recover its high-frequency component.

More precisely, we introduce the low and high-frequency projections, for some large $n$  to be fixed later on,
$$\P_n=\mathds{1}_{[0, n]}\big((AA^*)^{1/2}\big)\text{ and }\Q_n=I-\P_n.$$
For $U=U(t)$ a mild solution of \eqref{e:UCPabstraite} in $C^0([0, T], X^\sigma)$, let us consider the splitting
\begin{align*}
    U=\P_n U+\Q_nU=:U_{L}+U_{H}.
\end{align*}
Applying these projectors to the equation, noticing that $A$ and the projectors commute, we see that these components solve 
\begin{align}\label{e:lowhigh}
    \left\{\begin{array}{rl}
         \partial_t U_{L}(t)&=AU_{L}(t)+\Pc_n^{L} F(U_{L}+U_{H}),  \\
         \partial_t U_{H}(t)&=AU_{H}(t)+\Pc_n^{H} F(U_{L}+U_{H}), \\
\mathbf{C} U_{H}(t)     &=-   \mathbf{C} U_{L}(t).
    \end{array}\right.
\end{align}

We want to write $U_{H}=\mathcal{R}(U_{L})$ for an analytic reconstruction operator $\mathcal{R}$ and $n$ large. This will indeed be the case, see Theorem \ref{thmabstract} below. Assuming this for the moment, the equation by the low-frequency component in \eqref{e:lowhigh} writes
\begin{align*}
         \partial_t U_{L}(t)&=AU_{L}(t)+\Pc_n^{L} F(U_{L}+\mathcal{R}(U_{L})).
\end{align*}
This is now an autonomous equation for $U_{L}$. Moreover, since $A$ is a bounded operator and $\mathcal{R}$ is an analytic operator, this is an ODE in Banach space with an analytic nonlinearity. So $U_{L}$ is analytic in time. Therefore, it is also the case for $U_{H}=\mathcal{R}(U_{L})$ by composition of analytic functions. This gives the expected result once we have constructed $\mathcal{R}$.

The main task is therefore to build this reconstruction operator, as stated below.

\paragraph{Nonlinear reconstruction}
\begin{theorem}[Nonlinear reconstruction, Laurent-Loyola~\cite{LL:24}]
\label{thmabstract}Let $R_{0}>0$. There exist $n\in \N$ and a nonlinear Lipschitz (reconstruction) operator $\mathcal{R}$ 
\begin{align*}
\mathcal{R}:  \quad B_{R_{0}}(C^{0}([0,T],\Pc_n^{L}X^{\sigma}))\longrightarrow B_{R_{0}}(C^{0}([0,T],\Pc_n^{H}X^{\sigma}))
\end{align*}
 so that for any $U\in  B_{R_{0}}(C^{0}([0,T], X^{\sigma}))$ that satisfies
\begin{align*}
    \left\{\begin{array}{cl}
         \partial_t U=AU+F(U), &\textnormal{ on } [0,T], \\
         \mathbf{C}  U(t)=0, &\textnormal{ for }t\in [0,T],
    \end{array}\right.
\end{align*}
then, $\Pc_n^{H}  U=\mathcal{R}(\Pc_n^{L}U)$.

Moreover, $\mathcal{R}$ extends holomorphically in a small neighbourhood.
\end{theorem}
Note that this reconstruction operator is defined for any function, but is of interest only for the pathological functions satisfying the zero observation property \eqref{e:UCPabstraite}.

To prove Theorem \ref{thmabstract} and achieve the reconstruction, we need to inspect the equation satisfied by the high-frequency component of $\Pc_n^{L}U$, namely the last two lines of \eqref{e:lowhigh}:
\begin{align}
\label{e:high}
    \left\{\begin{array}{rl}
         \partial_t U_{H}(t)&=AU_{H}(t)+\Pc_n^{H} F(U_{L}+U_{H}), \\
\mathbf{C} U_{H}(t)     &=-   \mathbf{C} U_{L}(t).
    \end{array}\right.
\end{align}
We need to prove that for a fixed $ U_{L}$, regarded as a parameter or source term from now on, there is a unique solution to \eqref{e:high}. The first equation is a classical semilinear equation with a source term $ U_{L}$. So, if it were posed as an initial value problem, this would follow from the usual fixed point argument based on the Duhamel formula. However, in this setting, the low-frequency component is prescribed, whereas the initial value of the high-frequency component is not.

The idea is therefore to replace the usual initial value problem $U_{H}(0)=U_{H,0}$ by the observation condition $\mathbf{C} U_{H}=-\mathbf{C} U_{L}$ since the low frequency term $U_{L}$ is assumed to be known in this part of the problem. This suggests that given $U_{L}$, we could try to reconstruct $U_{H}$ by solving the corresponding nonlinear observability system. Indeed, in the linear case, according to \eqref{observabs}, the observability of the linear semigroup $t\in [0, T]\mapsto e^{tA}$ enables us to recover an initial condition $W_0$ solely in terms of an observation. Neglecting for the moment the source term given by the nonlinearity, this would allow us to reconstruct $W_0$ in terms of the observation of $W$, namely $\bC W=-\bC V$. Lemma \ref{lmCauchyobs} below provides a generalization of this reconstruction problem when source terms are present. It will take the role of the usual Duhamel formula but replacing the initial value problem by an observability problem.

\paragraph{Linear reconstruction}

\begin{lemma}[A Duhamel formula/observation problem]
\label{lmCauchyobs}
For any $H\in L^{1}([0,T],X^{\sigma})$ and $G\in L^{2}([0,T],X^{\sigma})$, there exists a unique $W\in C^{0}([0,T],X^{\sigma})$ solution of 
\begin{align*}
    \left\{\begin{array}{rl}
        \partial_t W(t)&=AW(t)+H,  \\
   \Pi  \mathbf{C} W    &=\Pi G. 
    \end{array}\right.
\end{align*}
\end{lemma}
Here $\Pi$ denotes the orthogonal projector onto the range of the observation operator $\mathbf{C}$. The necessity of introducing this projector can be easily seen by considering the homogeneous case $H=0$, which corresponds to reconstructing a free solution from its observation. It is clear then that not all observations are admissible. For instance, for a wave equation and an observation on $[0,T]\times \omega$, only those observations that are solutions of the free wave equation on $[0,T]\times \omega$ can be obtained from an observation of a free wave equation. Concerning the uniqueness, we are also led to the free wave equation and it is a direct consequence of the observability inequality. The proof of Lemma \ref{lmCauchyobs} relies mainly on functional analysis arguments, using the fact that the observability inequality~\eqref{observabs} ensures that the observation operator is injective and has closed range. Moreover, since we are interested only in the equation at high-frequency, we will require the following slight variant.
\begin{lemma}
\label{lmCauchyobsn}
For any $n\in\N$, $H\in L^{1}([0,T],X^{\sigma})$ and $G\in L^{2}([0,T],X^{\sigma})$, there exists a unique $W\in C^{0}([0,T],\Pc_n^{H}X^{\sigma})$ solution of 
\begin{align*}
    \left\{\begin{array}{rl}
        \partial_t W(t)&=AW(t)+\mathcal{P}_{n}^{H}H,  \\
   \Pi_n  \mathbf{C} W    &=\Pi_{n}G. 
    \end{array}\right.
\end{align*}
\end{lemma}
Here $\Pi_{n}$ denotes the orthogonal projector onto the range of the observation operator $\mathbf{C}$ restricted to $\mathcal{P}_{n}^{H}$.

\paragraph{End of the proof}

To complete the proof of Theorem \ref{thmabstract}, the idea is to carry out a fixed point argument for equation \eqref{e:high}, in a similar way as for a usual initial value problem, but using Lemma \ref{lmCauchyobsn} as a replacement for the Duhamel formula. Recall that $U_L$ appears as a source term and in the observation. It is supposed to be known in this reconstruction argument. 

In order to perform a fixed argument, some smallness is required. Usually, in the initial value problem, this is achieved by taking the time $T$ to be small or assuming small initial datum. In this context, neither option is available, since $T$ is fixed and the solution may be large. Instead, smallness is obtained by choosing the parameter parameter $n$ sufficiently large, corresponding to a high-frequency regime. In this framework, the third part \eqref{Fcomp} of Assumption \ref{assum}, namely, that $F$ is a nonlinear map from bounded sets of $X^\sigma$ into $X^{\sigma+\veps}$, provides the required smallness for large $n$. Therefore, at sufficiently high frequency $n$, we can treat the nonlinearity as a perturbation and successfully complete this reconstruction procedure. More precisely, this yields the existence of the nonlinear reconstruction map $\mc{R}$.

Finally, the analyticity of $\mc{R}$, which is the last claim of Theorem \ref{thmabstract}, is shown by checking the analyticity of all the terms involved. This relies on the assumption  that $F$ is analytic, that is part \eqref{Fanalytic} in Assumption \ref{assum}, and on standard results of the regularity of a fixed point with respect to parameters. As precised earlier this is actually essential to demonstrate that $t\in (0, T)\mapsto U_L(t)\in X^\sigma$ is analytic, and consequently that the same holds for $U_H$.
\section{Nonlinear plate equation}
\label{s:plates}

The previously described abstract framework allow us to obtain similar result for other PDEs, such as the nonlinear plate equation,
\begin{align}\label{eq:nlp}\tag{NLP}
    \left\{\begin{array}{rl}
        \partial_t^2 u+\Delta_g^2 u+\beta u+f(u)=0   &\ (t, x)\in [0, T]\times \M,\\
        u_{|_{\partial\M}}=\Delta u_{|_{\partial\M}}=0      &\ (t, x)\in [0, T]\times\partial\M,\\
        (u, \partial_t u)(0)=(u_0, u_1) &\ x\in \M,
    \end{array}\right.
\end{align}
where $f: \R\to \R$ is the nonlinearity. The following assumption will be made on the observation zone.

\begin{assump}{2}\label{assumOBS}
    The Schr\"odinger equation is observable in $L^2$ from $\omega$ in time $T_0>0$: there exists $C=C(T_0,\omega)>0$ such that for any $v_0\in L^2(\M)$ it holds
    \begin{align*}
        \norm{v_0}_{L^2(\M)}^2\leq C\int_0^{T_0} \norm{\mathds{1}_\omega e^{it\Delta_g}v_0}_{L^2(\M)}^2dt.
    \end{align*}
\end{assump}

Although it is an open question to find the optimal geometric condition on $\omega$ for which Assumption \ref{assumOBS} is verified, it is known to hold true in several (quite nontrivial) situations: $\omega$ satisfying the \eqref{assumGCC} on compact manifolds with or without boundary \cite{Leb92}, $\omega$ open touching the boundary in the two-dimensional euclidean disk \cite{2016:anantharaman-leautaud-macia:wigner-schrodinger-disk}, $\omega$ any open set in $\mathbb{T}^d$ \cite{K:92,AM:14}, $\omega$ any open set in hyperbolic surfaces \cite{DJN:22}, just to mention some. See \cite{BZ:04, AR:12} for further examples.
The unique continuation result is the following.
\begin{theorem}[Loyola \cite{Loy25:NLP}]
    Let $f$ be real analytic with $f(0)=0$. Let $0<T_0<T$ and $\omega\Subset\widetilde{\omega}$ are two nonempty open sets such that $\omega$ satisfies Assumption~\ref{assumOBS} in time $T_0$. If one solution $U=(u,\partial_t u)\in C^0([0, T], H^2\cap H_0^1\times L^2(\M))$ of \eqref{eq:nlp} satisfies $\partial_t u=0$ in $[0, T]\times \widetilde{\omega}$, then $\partial_t u=0$ in $[0, T]\times \M$ and $u$ is a stationary solution (solution of elliptic equation). If, moreover, the nonlinearity $f$ satisfies the defocusing assumption \eqref{defdefocusing}, then $u\equiv 0$.
\end{theorem}
As described in Section \ref{s:NLisL}, our strategy is based in the following:
\begin{itemize}
    \item The first ingredient is a recent unique continuation result for the linear plate operator due to Filippas-Laurent-Léautaud \cite{FLL:24}. It allows lower-order terms assuming Grevrey-$2$ regularity in time (in particular, analyticity fits in) and fairly general global geometric assumptions. This has to be compared with the Tataru-Robbiano-Zuily-Hörmander theorem, which notably does not apply to the plate operator \cite[Appendix B]{FLL:24}.
    \item The second ingredient is a propagation of analyticity in time for the nonlinear equation, obtained in Laurent and Loyola \cite{LL:24}. The nonlinear plate indeed fits in the abstract framework of Section \ref{s:abstract} and thus Theorem \ref{thmabstractanalyticintro} applies.
\end{itemize}

This unique continuation result is the main ingredient in \cite{Loy25:NLS} that allows us to obtain some semiglobal stabilization and controllability results for the nonlinear plate equation under more natural geometric conditions coming from Assumption~\ref{assumOBS}. Notably, most previously known results in the literature hold under strong assumptions of 'multiplier-type', see Section~\ref{s:carlemanhormander}. In this direction, we improve upon the current literature on unique continuation for nonlinear plates.

\section{Nonlinear Schr\"odinger equation}
\label{s:NLS}

From now on let us assume that $\M$ is a compact boundaryless Riemannian manifold of dimension $d$, and let us consider the nonlinear Schrödinger equation 
\begin{align}\label{I:eq:NLS}\tag{NLS}
        \left\{\begin{array}{cl}
            i\partial_tu+\Delta_g u=f(u)     & (t, x)\in (0, T)\times\M,  \\
            u(0)=u_0    &  x\in\M,
        \end{array}\right.
\end{align}
with nonlinearity $f: \mathbb{C}\to \mathbb{C}$. In a similar spirit as before, we have the following propagation of analyticity in time from a subdomain $(0, T)\times\omega$ where the \eqref{assumGCC} holds.

\begin{theorem}[Propagation of analyticity, Loyola \cite{Loy25:NLS}]
Let $d\in \N$ and $s>d/2$. Let $u\in C^0([0, T], H^{s}(\M))$ be a solution of \eqref{I:eq:NLS}. Assume :
    \begin{enumerate}
        \item $(\omega, T_0)$ satisfies the Geometric Control Condition~\eqref{assumGCC} for some $T_0>0$,
        \item $t\in (0, T)\mapsto \chi u(t,\cdot)\in H^{s}(\M)$ is analytic for any cutoff function $\chi\in C_c^\infty(\M)$ whose support is contained in $\omega$.
        \item $s\mapsto f(s)$ is real analytic with $f(0)=0$.
    \end{enumerate}
    Then $t\mapsto u(t, \cdot)$ is analytic from $(0, T)$ into $H^2(\M)$.
\end{theorem}

As described in Section \ref{s:UCPpartial}, we can obtain a unique continuation statement using the Tataru-Robbiano-Zuily-Hörmander result for the linear Schrödinger equation with potential; see \cite[Theorem 6.4]{LL:19} for a quantitative statement of such result. In the $H^1$-subcritical case, we assume that $(\M, g)$ can be any of the following manifolds:
\begin{itemize}
    \item a compact boundaryless surface;
    \item $\mathbb{T}^3$ or the irrational torus $\R^3/(\theta_1 \Z \times \theta_2 \Z \times \theta_3\Z)$ with $\theta_i \in \R$;
    \item $S^3$ or $S^2\times S^1$.
\end{itemize}
Furthermore, we consider $f(u)=P'(|u|^2)u$ to be of polynomial type, where:
\begin{enumerate}
        \item\label{I:NLS:reg:2.1} if $d=2$ then $P$ is a polynomial function with real coefficients, satisfying $P(0)=0$ and the defocusing assumption $P'(r)\xrightarrow[r\to+\infty]{}+\infty$;

        \item\label{I:NLS:reg:2.2} if $d=3$, then $P'(r)=\al r+\beta$ with $\al>0$, $\beta\geq 0$, corresponding to the cubic nonlinearity.
\end{enumerate}

The unique continuation result is the following.

\begin{theorem}[Unique continuation, Loyola \cite{Loy25:NLS}]
    Let $(\M, g)$ be any of the manifolds described above and let $f$ be as in \ref{I:NLS:reg:2.1}-\ref{I:NLS:reg:2.2} according to the dimension of the manifold. Assume that there exists some $T_0>0$ such that $(\omega, T_0)$ satisfies the Geometric Control Condition. If one solution $u$ of \eqref{I:eq:NLS}, belonging to a suitable functional space $\mc{X}_d\subset C([0, T], H^1(\M))$, satisfies $\partial_tu=0$ in $(0, T)\times\omega$, then $u=0$ in $(0, T)\times \M$.
\end{theorem}

Here $\mc{X}_d$ denotes a suitable functional space inherited from the wellposedness framework in each case. Namely, for $d=2$, $\mc{X}_2$ ensures finite Strichartz norm and for $d=3$, $\mc{X}_3$ corresponds to the Bourgain space $X_T^{1, b}$ with $b\in (1/2, 1]$. This answers in the affirmative an open question posed by Dehman-Gérard-Lebeau \cite{DGL06} in the nonlinear case.

Notably, the compactness assumption on the nonlinearity, namely, part \eqref{Fcomp} of Assumption~\ref{assum}, does not hold for \eqref{I:eq:NLS}. We therefore need a variant of Theorem \ref{e:UCPabstraite} that fits this framework. Below we briefly explain how to drop the compactness assumption on the nonlinearity, stressing the main issues in the proof.

\subsection{On the nonlinear reconstruction}

The set of assumptions can be roughly stated as follows. 

\begin{assump}{3}\label{assumNLS}
    \begin{enumerate}
        \item $A$ is skew-adjoint with compact resolvent, giving ''Sobolev spaces'' $X^\sigma$.
        \item\label{FanalyticLip} $F: X^\sigma \to X^\sigma$ is analytic, and both $F$ and $DF$ are Lipschitz on bounded subsets.
        \item\label{unifOBS} Uniform observability: for $\bC\in\mc{L}(X^\sigma)$, a large $n$ and suitable inputs $V\in\mathbb{V}$,
        \begin{align*}
            \left\{\begin{array}{c}
                \partial_t W=(A+\Q_n DF(V))W, \\
                W(0)=W_0\in \Q_n X^\sigma, 
            \end{array}\right.\ \Longrightarrow\ \norm{W_0}_{X^\sigma}^2\leq \mathfrak{C}_{\text{obs}}^2\int_0^T \norm{\bC W(t)}_{X^\sigma}^2dt,
        \end{align*}
        with $\mathbb{V}$ a technical set which enjoys some compactness properties.
    \end{enumerate}
\end{assump}

In this setup, the following variation of the abstract Theorem \ref{thmabstractanalyticintro} holds.

\begin{theorem}
Let $T^{*}>T$ and $\mc{K}$ be a "suitable" compact set of $C^0([0, T^*], X^\sigma)$. Then any solution $U\in \mc{K}$ of
    \begin{align}\label{NLS:eq:UCabstT*intro}
        \left\{\begin{array}{cr}
             \partial_t U=AU+F(U)& \textnormal{ on } [0,T^{*}], \\
             \bC U(t)=0 &\textnormal{ for }t\in [0,T^{*}],
        \end{array}\right.
    \end{align}
    is real analytic in $t$ in $(0,T^*)$ with value in $X^{\sigma}$.
\end{theorem}

The key issue lies in how to perform a similar reconstruction procedure as in Section \ref{s:scheme}. For a solution $U$ to \eqref{NLS:eq:UCabstT*intro}. We split in low and high frequencies $U=U_L+U_H$ and by studying the equation satisfied by the high-frequency component, we aim to find a suitable analytic reconstruction operator $\mc{R}$ so that $U_H=\mc{R}(U_L)$. To face the lack of compactness of the nonlinearity $F$, we take into account the linear variation of the high-frequency component $U_H$ along $\Q_nDF(U_L)$ when $U$ is a solution living on a compact set. That is, we are led to study the high-frequency observability problem:
\begin{align*}
    \left\{\begin{array}{l}
         \partial_t U_H(t)=\big(A+\Q_nDF(U_L)\big)U_H+\Q_n\mc{H}(U_L, U_H), \\
         \bC U_H(t)=-\bC U_L(t).
    \end{array}\right.
\end{align*}
where $\mc{H}(U_L, U_H)=F(U_L+U_H)-DF(U_L)U_H$. We stress three key issues:

\begin{enumerate}
    \item Part \eqref{unifOBS} of Assumption \ref{assumNLS} allows us to perform the reconstruction for the linearized system at high-frequency, similar to Lemma \ref{lmCauchyobsn} and considering the projector accordingly. This is why we enforce the observability estimate to be uniform with respect to the frequency threshold and to low-frequency parameter as well. 
    \item At high-frequency $n$, the Lipschitz control provided by $F$ (part \eqref{FanalyticLip} of Assumption \ref{assumNLS}) and that the solution itself enjoy some compactness properties, allows us to close the fixed point argument. This guarantees the existence of a reconstruction operator $\mc{R}$.
    \item The analyticity of $\mc{R}$ is delicate as it depends in a (quite)nonlinear way with respect to the low-frequency input: not only through the nonlinearity $F$, but also through the linearization at high-frequency and the corresponding projector of the linearized observed solutions.
\end{enumerate}

Once we are able to write $U_H=\mc{R}(U_L)$ with $\mc{R}$ analytic, the proof follows similarly as described in Section \ref{s:scheme}.


\end{document}